\documentstyle[11pt]{article}
\newcommand{\be}{\begin{equation}}
\newcommand{\ee}{\end{equation}}
\newcommand{\bea}{\begin{eqnarray}}
\newcommand{\eea}{\end{eqnarray}}
\newcommand{\ra}{\rightarrow}

\newcommand{\Lra}{\Longleftrightarrow}
\newcommand{\lra}{\longleftrightarrow}
\newcommand{\Ua}{\Uparrow}

\def\Journal#1#2#3#4{{#1} {\bf #2}, #3 (#4)}

\def\PLB{{\em Phys. Lett.} B}

\def\CMP{\em Commun. Math. Phys.}
\def\CTP{\em Commun.Theor. Phys.}
\def\LMP{\em Lett. Math. Phys.}
\def\SJNP{\em Sov. J. Nucl. Phys.}
\def\JMP{\em J. Math. Phys.}     

\def\JPA{{\em J. Phys.} A}
\def\JPG{{\em J. Phys.} G}
\def\PPP{\em Prog.Part.Nucl.Phys.}
\def\RPP{\em Rep.Prog.Phys.}
 
\begin{document}
\large
\baselineskip=18pt

\thispagestyle{empty}

\hfill {\small DSF-44/01}

\vspace{10mm}

\begin{center}
{\Large \bf A $q$-analogue of the embedding chain $U(6) \supset G
\supset S0(3)$} \\[6mm]
 A. Sciarrino\\[3mm]
Universit\`a di Napoli ``Federico II'' \\
Dipartimento di Scienze Fisiche \\
and \\
I.N.F.N. - Sezione di Napoli \\
I-80126 Napoli - Italy
\end{center}

\bigskip

\bigskip

\bigskip

\bigskip

 \centerline{\bf Abstract}

\begin{quote}
 A $q$-analogue of the embedding chains of the Arima-Iachello model is
 proposed. The generators of the deformed $U(6)$ subalgebras are 
 written in terms of the generators of $gl_{q}(6)$, using $q$-bosons.
\end{quote}

 \vfill
 {\small \bf
\begin{tabular}{l}
~$\,$Postal adress:  Complesso Universitario di Monte S. Angelo \\ 
~$\,$Via Cintia - I-80126 Napoli (Italy)\\
 ~$\,$E-mail: sciarrino@na.infn.it 
 \end{tabular} }

\pagebreak
\thispagestyle{empty}
\mbox{}
\newpage
\setcounter{page}{1}
\section{Introduction}

 Since their introduction \cite{J}, \cite{D}, the quantum algebras $G_q$ or 
 $U_q (G)$, i.e. the $q$-deformed universal 
enveloping algebra of a semi-simple Lie algebra $G$ have been  a topic of active
 research both in physics and mathematics. The underlying idea in some 
 of their applications  is to use a q-deformed algebra
instead of a Lie algebra to realize a {\em generalized dynamical symmetry}.
For a review of the applications and methods of the dynamical 
 or spectrum generating algebras in physics, see 
\cite{Gruber}, and, in nuclear physics, \cite{VanI}.
The key idea of dynamical symmetry scheme is to write the hamiltonian 
of  a physical system as a sum of invariants, usually second order 
Casimir $C$, with constants to be determined by experimental data, of the 
embedding chains of algebras of the type:
\be
G \supset L \supset \ldots \supset SO(3)    \label{eq:0}
\ee
\be
{\cal H} = C(G) \, + \, C(L) \,  \ldots \, + C(SO(3))  
\ee
\noindent where $SO(3)$ describes the angular momentum and, usually, 
the Casimir operators are  written  using Jordan-Schwinger like 
realization of the  algebra $G$ by means of bosonic creation-annihilation 
operators. The idea of dynamical symmetry has countless applications in 
molecular, atomic, nuclear, hadronic and chemical physics.
The most simple example is the rigid rotator where the Hamiltonian is
written as the Casimir operator of $SU(2)$
\be
C = k \, J(J \, + \, 1)  \label{eq:rr}
\ee
 The energy spectrum of eq.(\ref{eq:rr}) is of the form
\be
E_{j} = k \, j(j + 1) \label{eq:srr}
\ee
The replacement in eq.(\ref{eq:rr}) of the Casimir of $SU(2)$ by the
Gasimir operator of $sl_{q}(2)$ \cite{RRS} provides the first example of 
application of deformed algebra  as  dynamical symmetry (see below for 
 notation)
\be
C = K \, [J]_{q}[J \, + \, 1]_{q}  \label{eq:drr}
\ee
Now the energy spectrum will depend on a parameter $ q = \exp{i \tau}$
and we have 
\be
E_{j}^{q} = K \, [j]_{q}[(j + 1)]_{q} = K \, \frac{\sin(|\tau| 
j)\sin(|\tau(j + 1))}{\sin^{2}(|\tau|j)} \label{eq:sdrr}
\ee
 Eq.(\ref{eq:sdrr}) fits  experimental data for several deformed nuclei, \cite{RRS}, 
 better than  eq.(\ref{eq:srr}) ($|\tau| = 0$). The results of this simple model 
 suggests that it may be worthwhile to further investigate the idea of 
 generalized dynamical symmetry based on $q$-algebras. Indeed 
 in this last decade many applications, mainly, but not uniquely, in molecular and 
 nuclear physics have been investigated. For an excellent review of the subject with 
 an exhaustive list of references see \cite{BD}.  
The simplest, not trivial, embedding chain is the so called Elliot 
model: 
\be
 SU(3) \supset SO(3)  \label{eq:Em}
 \ee 
The deformation of the simple embedding chain of eq.(\ref{eq:Em})  is 
not trivial. Indeed in \cite{Van} it has been shown that the generators 
of $ so_q(3)$ can be expressed by means of the generators of  
$gl_q(3)$, not of $sl_q(3)$, iff the algebraic relations are restricted to the 
symmetric representations.  Moreover the coproduct of $gl_q(3)$ 
{\bf does not induce} the standard coproduct on  $ so_q(3)$.
It is useful to emphasize that the definition of the coproduct is essential 
to define the tensor product of spaces. The $q$-analogue of the chain 
eq.(\ref{eq:Em}) has been widely studied. In  \cite{Scia1} the author 
has proposed a possible solution, but the problem has been tackled 
from several points of view, see \cite{BD} for references to the 
different solutions and for physical applications of the $q$-analogue 
of the embedding chain eq.(\ref{eq:Em}).
 Therefore it is clear that an essential step to carry forward the program of 
 application of q-algebras as { \em generalized dynamical symmetry}, beyond the
  simple models  above discussed, is to dispose of a formalism which 
allows to build up analogous chains of eq.(\ref{eq:0})  replacing the Lie 
algebras by the deformed ones. Of course, as we are no
more dealing with Lie algebras, the term {\bf embedding} has to be intended
in the loose sense that the generators of the embedded deformed subalgebra are
expressed in terms of the generators of the algebra while the
Hopf structure can be or inherited from that of the embedding algebra
or imposed on the generators of the embedded algebra. 
The root of the problem,as it has been discussed in \cite{Scia2}, lies on the 
fact that $G_q$ are well defined only
in the Cartan-Chevalley basis and this basis is not suitable to discuss 
embedding of subalgebras except the regular ones. The classification 
of so called singular subalgebras of Lie algebras has been started by 
Dynkin and we refer to the clear paper of Gruber and Lorente \cite{GL},  where the
embedding matrices are explicitly computed for low rank algebras.
In particular in \cite{Scia2} it has been shown that, in the case where the rank of $L$,
maximal singular algebra of $G$, is equal to the rank of $G$ minus one, it
is possible, using realization of $G_q$ in terms of $q$-bosons and/or in
terms of the so called $q$-fermions, to write the Cartan-Chevalley generators
of $L_q$ in terms of the generators of $G_q$. Let us remark that this result is
not at all a priori obvious due to the non linear structure of $G_q$.
It has also been discussed
what kind of deformed $G$ is obtained if the standard coproduct is imposed
on the generators of $L_q$ in the standard way instead of being derived from
that of $G_q$. 
 The aim of this paper is to focus on the embedding 
chains appearing in the  Arima-Iachello model \cite{AI}, \cite{IA}
and to discuss in which sense one can write  analogous  chains of $q$-algebras.  
This very successfull model is based on the following three embedding chains
\be
  SU(6) \supset \left \{ \begin{array}{ll}
  SU(5) \supset SO(5) \supset SO(3) & \mbox{(vibrational)} \\
   SU(3) \supset SO(3) & \mbox{(rotational)} \\
    S0(6) \supset SO(5) \supset SO(3) & \mbox{($\gamma$-unstable)}
    \label{eq:AI}
    \end{array} \right.
    \ee 
A partial answer to the question of finding a $q$-analogue of the above 
embedding chains has been given in \cite{NS}, where it 
has been shown that, using $q$-bosons realization, the deformed maximal $SU(6)$
subalgebras, i.e. $sl_q(3)$,
$so_q(6)$, can be written in terms of the generators of $gl_q(6)$ and that
this procedure can be extended also to the  deformation of $SO(5)$, maximal subalgebra
of $ SU(5) \subset SU(6) $.
It should be mentioned that an attempt of writing a deformed version of the
Arima-Iachello model has been carried out in \cite{WY}, using the notion of
complementary subalgebras introduced about thirty years ago by 
Moshinsky-Quesne \cite{MQ}. Two subalgebra $L_{1}$ and $L_{2}$ of an algebra
$G$ are complementary in one definite irrep. of $G$, if there is a 
correspondence one-to-one between the irreps. of  $L_{1}$ and $L_{2}$ 
contained in the considered irrep. of $G$. Using this notion  in \cite{WY} 
 a hamiltonian has been written in terms of the second order Casimir of
 $su_{q}^{sd}(1,1)$,  $su_{q}^{d}(1,1)$ and  $su_{q}(2)$, where $s,d$ are boson 
 operator and the second  $su(1,1)$  is contained in the first one.
  This hamiltonian, in the limit $q \ra 1$, 
 for particular value of the coefficients of the Casimir 
 operators, tends to a hamiltonian in the $\gamma$-unstable chain.
 Althoug this approach   is interesting, it should be pointed out that
 the content and the embedding of subalgebras is not well defined. Rather 
 than a deformation of the embedding chain eq.(\ref{eq:AI}), it looks like 
 an interesting dynamical model based on $q$-deformed algebras. 
  Our aim is to build up the whole algebraic construction of a 
  $q$-analogue of the chains of eq.(\ref{eq:AI}). From the
  construction it should be possible to write hamiltonians which in the
  limit $q \ra 1$ tend to hamiltonians of the undeformed model,  for any
  value of the coefficients of the Casimir operators. Indeed we shall show
  that, replacing $U(6)$ by  $gl_{q}(6)$,  we can write the generators of any
  $U(6)$ subalgebra in terms of the generators of $gl_{q}(6)$.
To make the paper self-contained,  in Sec. 2 we  briefly recall
 the formulas we shall use in the following. In Sec. 3 we write explictly the 
$q$-analogue of the embedding chains eq.(\ref{eq:AI}). In Sec. 4
we summarize and discuss our results.  
 
 \section{Reminder }

 We recall, also to fix the notation, the definition
of deformed Lie algebra $G_q$ in the Cartan-Chevalley basis, of $sl_{q}(2)$ $q$-tensor 
operators, the definition
of q-bosons, which we shall use to write explicit realization of the q-algebras,
 the  Curtright-Zachos deforming functional, connecting $SU(2)$ and 
$sl_{q}(2)$ and the deforming map between $Sp(4) \equiv SO(5)$ and 
$sp_{q}(4) \equiv so_{q}(5)$. 

\subsection{Deformed algebras}

Let us recall the definition of $G_q$ associated with a simple Lie algebra
$G$ of rank $r$ defined by the Cartan matrix $ (a_{ij}) $ in the Chevalley
basis. $G_q$ is 
generated by $3r$ elements $e_i^{\pm}$  and $h_i$ which satisfy 
$(i,\, j=1= 1, \ldots , r)$

\be \begin{array}{c}
~[e_i^+,\, e_j^-]  =  \delta_{ij} [h_i]_{q_i} ~~~~~~~~~~~~~~~~ [h_i, h_j] = 0 \\
~~ \\
~[h_i, \, e_j^+]  =  a_{ij} e_j^+ ~~~~~~~~~~~~~~~~~~ [h_i, e_j^-] 
= -a_{ij} e_j^- 
\end{array}
\ee
where 
\be 
~[x]_q = \frac{q^x - q^{-x}} {q - q^{-1}}  
\ee
and $q_i = q^{d_i}$, $d_i$ being non-zero integers with greatest 
common divisor equal to one such that $d_i \, a_{ij} = d_j\, a_{ji}$.
For simple laced algebras $d_i = 1$ while for $so_{q}(2n+1)$ ($sp_{q}(2n)$)
$d_i = 2 \;\; (1)$, $i \neq n$, $d_n = 1 \;\; (2)$.  
Further the generators have to satisfy the Serre relations:
\be
\sum_{0 \leq n \leq 1 - a_{ij}} ~ (-1)^n \left[ \begin{array}{c}
1 - a_{ij} \\ n \end{array} \right]_{q_i} \, (e_i^+)^{1 - a_{ij} - n} \, 
e_j^+ (e_i^+)^n =0 
\ee
where
\be
\left[ \begin{array}{c} 
m \\ n \end{array} \right]_q = \frac{[m]_q!}{[m - n]_q! [n]_q}
\ee
\begin{center}
$[n]_q! = [1]_q \: [2]_q \ldots [n]_q$
\end{center}
Analogous equations hold replacing $e_i^+$ by $e_i^-$. In the following
we assume  $~h_i~$ = $~(h_i)^{\dag}~$ and the deformation parameter $q$ to be 
different from the roots of the unity.
The algebra $G_q$ is endowed with a Hopf algebra structure, i.e. on $G_q$  
the action of the coproduct $\Delta$, antipode $S$ and co-unit $\varepsilon$ 
 is defined. This extremely relevant aspect will not discussed here. 
 We recall only the definition of the coproduct  which we  
 shall briefly refer to in the following.
\be
\begin{array}{c}
\Delta (h_i) = h_i \otimes {\bf 1} + {\bf 1} \otimes h_i \\
~~ \\
\Delta (e_i^{\pm}) = e_i^{\pm} \otimes q_i^{h_i/2} + 
q_i^{-h_i/2} \otimes e_i^{\pm} 
\end{array}  \label{eq:co}
\ee
\noindent As the coproduct in a Hopf algebra satisfies ($ g_i, g_j \in G_q $)
\be
\Delta(g_i g_j) = \Delta(g_i) \Delta(g_j)
\ee	
\noindent it is essential to define which elements \{$g_i$\} are the ``basis" 
of $G_q$.

Let us recall the definition of $gl(n)_{q}$,  ($i,j = 1,2, \dots, n-1, k = 
1,2, \dots, n-1,n$):
\be \begin{array}{c}
~[e_i^+,\, e_j^-]  =  \delta_{ij} [n_{i} - n_{i+1}]_{q} ~~~~~~~~~~~~~~~~ 
[n_i, n_k] = 0 \\
~~ \\
~[n_k, \, e_j^{\pm}]  = \pm ( \delta_{k,j} - \delta_{k-1,j} \,)   
\end{array}
\ee
The Serre relations are computed using $ a_{k,j} = -(\delta_{k-1,j} + 
\delta_{k,j-1} \,)$. So $gl(n)_{q}$ can be considered as $sl(n)_{q} 
\oplus n_{n}$  with $h_{i} = n_{i} - n_{i+1}$. 

\subsection{ $q$-Bosons}

Let us recall the definition of 
Biedenharn-MacFarlane $q$-bosons  \cite{B}, \cite{M} which we 
denote by  $~ b_i^+, b_i ~$ 
\be
b_i b_j^+ - q^{ \delta_{ij}} \, b_j^+ b_i = \delta_{ij} q^{- N_i} 
\ee
\be
~[N_i, \, b_j^+] = \delta_{ij} b_j^+ ~~~~~ [N_i, \, b_j] = - 
\delta_{ij} b_j ~~~~~ [N_i,\, N_j] = 0 
\ee	
It is useful to keep in mind the following identities:
\be
b_i^+ b_i = \frac{q^{N_i} - q^{-N_i}} {q - q^{-1}} ~~~~~~~~
b_i b_j^+ = \frac{q^{N_i+1} - q^{-N_i-1}} {q - q^{-1}}
\ee
\be
 b_i^+ b_k  = [b_i^+ b_j, \, b_j^+ b_k]_{q} \,  q^{N_j} \label{eq:qc}
\ee 
where  the  $q$-commutator is defined as
\be
[A, \,B ]_{q} = A \, B \; - \; q B \, A
\ee
 The explicit construction of $q$-bosons in terms of non-deformed 
standard bosonic oscillators ($\tilde{b}_i^+, \tilde{b}_i$)  is 
\cite{S}
\be
b_i^+ = \sqrt{\frac{[N_{i}]_{q}}{N_{i}}} \, \tilde{b}_i^+ ~~~~~~~~
b_i  = \tilde{b}_i \, \sqrt{\frac{[N_{i}]_{q}}{N_{i}}}  \label{eq:song}
\ee 
From the above equation, as it is well known that the deformed algebras
$su_{q}(n)$ and $sp_{q}(2n)$ can be written in terms of bilinears of $q$-bosons,
\cite{Ha}, it follows that one can write the following correspondences
$su_{q}(n) \subset gl_{q}(n) \, \Lra \, gl(n) \supset sl(n)$ and
$sp_{q}(2n)  \, \Lra \,  sp(2n)$, if  each generator of algebras and $q$-algebras 
is written as a bilinear of bosons and, respectively, of $q$-bosons.  Remark that, if
$\tilde{b}_i^+$ is the adjoint of $\tilde{b}_i$, then $b_i^+$ is the 
adjoint of $ b_i$ iff $q$ is real or $q = \exp{i \tau}$, $\tau$ real. 

\subsection{$q$-bosons realizations of $sl_{q}(2)$ and $so_{q}(3)$}

In order to clarify what we mean by $sl_{q}(2)$ and $so_{q}(3)$,
 let us write explicitly the q-boson realization of $sl_{q}(2)$ 
\be
J_{+} =  b_1^+ b_2 \;\;\;\;\;\;  J_{-} = b_2^+ b_1 
\;\;\;\;\;\;  2 J_{0} =  N_1  - N_2  
\ee
the states of the irreducible representation ($j, m$)  in the 
corresponding Fock space are
\be
\psi_{jm} = \, \frac{( b_1^+)^{j+m} \, (b_2^+)^{j-m}}{\sqrt{[j + 
m]_{q}![j- m ]_{q}!}} \; \psi_{0}
\ee 
and of $so_{q}(3)$
\bea
L_{+} &  =  & q^{N_{-1}} \, q^{-N_{0}/2} \, \sqrt{q^{N_{1}} + q^{-N_{1}}} \, 
b_1^+ b_0 \nonumber \\ 
& +  &  b_0^+ b_{-1} \,  q^{N_{1}} \, q^{-N_{0}/2} \, 
\sqrt{q^{N_{-1}} + q^{-N_{-1}}}   \nonumber \\ 
  L_{-} &  =  &  b_0^+ b_{1} \, q^{N_{-1}} \, q^{-N_{0}/2} \, 
  \sqrt{q^{N_{1}} + q^{-N_{1}}} \nonumber \\Ê 
  & +  & \,  q^{N_{1}} \, q^{-N_{0}/2} \, 
\sqrt{q^{N_{-1}} + q^{-N_{-1}}} \,  b_{-1}^+ b_{0}   \nonumber \\Ê 
    L_{0} & = & N_1  - N_{-1}  
\eea
the states of the odd dimensional irreducible representation ($L, m = n_{1} - n_{-1}$)
  ($ L = max\{ n_{1}\} = n_{1} + n_{0}  + n_{-1}$)  in the corresponding 
Fock space are linear combinations of 
\be
\psi_{n_{1},n_{0},n_{-1}} = \, \frac{( b_1^+)^{n_{1}} \, 
(b_0^+)^{n_{0}} \, ( b_{-1}^+)^{n_{-1}}}{\sqrt{[n_{1}]_{q}![n_{0}]_{q}!
[n_{-1}]_{q}!}} \; \psi_{0}
\ee 

\subsection {\bf q-Tensor operator}

From the formula of coproduct eq.(\ref{eq:co}) Biedenharm and Tarlini 
\cite{BT}, see also \cite{STK} and \cite{RS}, have derived the general structure
of q-tensor operators for $sl_{q}(2)$. 
\bea
   \left [  J_{\pm}, \, T_{m}^{k}(q) \right ]_{q^{-m}} \, q^{-J_{0}}  
 & = & \sqrt{[k \mp m]_{q}[k \pm m + 1]_{q}} \,  T_{m \pm 1}^{k}(q)   \nonumber \\  
 \left [ J_{0}, \, T_{m}^{k}(q) \right ] & = & m \, T_{m}^{k}(q)  
\label{eq:tensor}
\eea
The Wigner-Eckart theorem for $sl_{q}(2)$ reads:
\bea
<JM| T_{m}^{k}(q)|j_{1}m_{1}>  & = & (-1)^{2k} \, [2J + 1]_{q}^{-1/2} \, 
<JM|km,j_{1}m_{1}>_{q}  \nonumber \\  & \times & <J || T^{k}||j_{1}>  
\eea
where  $<JM|km, j_{1}m_{1}>_{q}$ is the q-Clebsch-Gordan coefficient, see 
\cite{STK}

\subsection {\bf   Deforming map between $sl(2)$ and $sl_{q}(2)$}

In \cite{CZ} an invertible deforming functional ${\cal Q_{\pm}}$ has been 
introduced which allows to  relate $sl(2) \Lra sl_{q}(2)$.
Denoting by small (resp. capital) letter $j_{\pm, 0}$ ($J_{\pm, 0}$) the generator
 of $sl(2)$  ($sl_{q}(2)$) it is possible to write ($q$ real)
 \be
J_{+} = {\cal Q}_{+} (j_{\pm},j_{0}) \, j_{+}  \;\;\;\;\;\; 
J_{-} = {\cal Q}_{-} (j_{\pm},j_{0}) \, j_{-}  \;\;\;\;\;\; 
J_{0} =  j_{0}  \label{eq:CZ}
\ee 
where
\be
{\cal Q}_{+}  = \sqrt{ \frac{ Ê\{[J_{0} + \mbox{\bf J}]_{q} 
[J_{0} - \mbox{\bf J} - 1]_{q}  \} }
{ Ê\{(j_{0} +  \mbox{\bf j}) (j_{0} -  \mbox{\bf j} - 1)  \} } }
\ee
(${\cal Q}_{-} = {\cal Q}_{+}^{\dag}$) and the operator {\bf j} ({\cal J}) is defined
 by the Casimir operator of $sl(2)$ ($sl_{q}(2)$)  
\be
C =  \mbox{\bf j} \, ( \mbox{\bf j}+ 1) \;\;\;\;\;\;\;\;\;\;\;\;  
\left ( C_{q} = [ \mbox{\bf J}]_{q} \, [ \mbox{\bf J} + 1]_{q} \right )
\ee

\subsection {\bf  Deforming map between $SO(5)$ and $so_{q}(5)$}

In \cite{Scia3}, in the space of the symmetric irreps.,
invertible deforming maps have been derived which allow to express
 $sl_{q}(n)$ and $sp_{q}(2n)$ respectively in terms of $U(n)$ and
 $Sp(2n)$. As irreducible representations in the Fock space of bosons or
 $q$-bosons are symmetric and $SO(5) \equiv  Sp(4)$, we report here
 explicitly the map $sp_{q}(4)\; \Lra \; Sp(4)  $ which we shall use in 
 the following.  Denoting  by small (capital) letters the generators of
 deformed (undeformed) algebra, we have
 \bea
& e^{+}_{1} = E^{+}_{1} \, \sqrt{\frac{[H_{1} \, + \, H_{2} \, + \, 1]_{q} \, 
[H_{2}]_{q}}{(H_{1} \, + \, H_{2} \, + \, 1) \, H_{2}}}    \nonumber \\
& e^{+}_{2} = \frac{2}{q + q ^{-1}} \, E^{+}_{2} \, \sqrt{\frac{[( H_{2}  +  1]_{q}  \,
[- H_{2}  -  2)]_{q}}{( H_{2}  +  1) \,(- H_{2}  -  2)}}  \label{eq:map}  
\eea
\be
 h_{k} =  H_{k} \;\;\;\;\;\;\;\;\;\;\;\;  k = 1, 2
 \ee
 
\section {\bf $q$-Embedding }
 
In this Section we  discuss in detail the meaning of the $q$-embedding 
chain, $L$ being a maximal singular subalgebra of $G$
\be
G_q \supset L_q  \supset so_q(3)  \label{eq:defsub}
\ee 
 for the case where $ G_q = gl_q(6)$.  .
In the following we  denote by small (resp. capital) letter $e^{\pm}, h$ 
($E^{\pm}, H$) the generators of $gl_{q}(6)$ ($L_{q}$) and by $L_{\pm,0}$
the generators of $so_q(3)$. With a hat on $E^{\pm}, H$ we denote, 
when required, the generators of a maximal subalgebra of $L$. Let us 
recall once again that eq.(\ref {eq:defsub}) holds if the generators of 
 $ L_q $ can be expressed in terms of the generators of $G_q $, at 
 least in a particular realization of $G_q$,  in the present work 
 using a $q$-boson realization.
 In the following we assume $q$ real, therefore any generator $X^{-}$ is 
 the adjoint of $X^{+}$ and the generators of the Cartan subalgebra 
 are self-adjoint. We shall comment in the Conclusions on the more 
 general case.  
  
  Let us recall the $q$-boson realization of $sl_{q}(6) \subset gl_{q}(6)$  
 ($ i = 1, 2, 3, 4, 5$)
 \be
e_{i}^{+}   =   b_i^+ b_{i+1}   \;\;\;\;\;\;\;\;\;\;\;\; 
  h_{i}  =  N_{i}  - N_{i+1}   \label{eq:real}
  \ee
  To get $gl_{q}(6)$ one has to add to the previous generators
  \be
   h_{0}  = \sum_{i=1}^{5} \,  h_{i}  = \sum_{j=1}^{6} \,  N_{j}
   \ee   

In the following we write explicitly the $q$-analogue of the embedding chains of 
the Arima-Iachello mode. The notation is self-explanatory. We  use a vertical arrow 
in the  equations to point out what the embedding chains tend to in the limit 
$q \ra 0$. 
In the following, to save space in some equations we shall  denote by 
$e^{+}_{i}$ the generators of $sl_{q}(6)$, whose content in $q$-bosons is
given in eq.(\ref{eq:real}). Let us remind also that the Cartan generators
of the deformed and undeformed algebras are the same.

\subsection{ \bf $q$-analogue of the vibrational embedding chain}  

\bea
gl_{q}(6) \supset gl_{q}(5)  \supset  so_{q}(5) & \nonumber \\
& \supset^{\Ua_{q \ra 1}}  so_{q}(3)  
\label{eq:dv}
\eea
Clearly the generators of $gl_{q}(5)$ are  obtained from those of 
$gl_{q}(6) $ neglecting  $e_{5}^{\pm}, h_{5}$ and $ N_{6}$   
and the $q$-boson realization of  $so_{q}(5) \subset gl_{q}(5) \subset gl_{q}(6)$
is \cite{NS}
\bea
E_{1}^{\dag}  &  =  & \left \{ \sqrt{q^{N_{1}} + q^{-N_{1}}}Ê\,  
b_1^+ b_2 \, \sqrt{q^{N_{2}} + q^{-N_{2}}} \, q^{-(N_{4} - N_{5})} 
\right. \nonumber \\  
 &  +  &  \left. \sqrt{q^{N_{4}} + q^{-N_{4}}}Ê\,  
b_4^+ b_5 \, \sqrt{q^{N_{5}} + q^{-N_{5}}}  q^{(N_{1} - N_{2})} 
\right \} (q + q^{-1})^{-1} \nonumber \\Ê 
E_{2}^{\dag}  &  =  &  q^{N_{4} - N_{3}/2} \,  \sqrt{q^{N_{2}} + q^{-N_{2}}}Ê\,  
b_2^+ b_3   
\, +  \, b_3^+ b_4 \, q^{N_{2} - N_{3}/2} \,\sqrt{q^{N_{4}} + q^{-N_{4}}}  
\nonumber \\Ê   
 H_{1} & = & N_{1}  - N_{2} + N_{4}  - N_{5} \nonumber \\  
  H_{2} & = & 2 \, (N_{2}   - N_{4})    
\eea
the generators of $so_{q}(3)$ can be written in terms of the 
generators of $gl_{q}(5)$
\bea
L_{+} & = & {\cal Q_{+}}  \left \{ 2 \left [ 
\sqrt{\frac{N_{1}}{[N_{1}]_{q}}} e^{+}_{1} \sqrt{\frac{N_{2}}{[N_{2}]_{q}}}
\, + \, \sqrt{\frac{N_{4}}{[N_{4}]_{q}}} e^{+}_{4} 
\sqrt{\frac{N_{5}}{[N_{5}]_{q}}} \, \right ]  \right. \nonumber \\  
& + &  \sqrt{6} \left.  \left [ 
\sqrt{\frac{N_{2}}{[N_{2}]_{q}}} e^{+}_{2} \sqrt{\frac{N_{3}}{[N_{3}]_{q}}}
\, + \, \sqrt{\frac{N_{3}}{[N_{3}]_{q}}} e^{+}_{3} 
\sqrt{\frac{N_{4}}{[N_{4}]_{q}}} \, \right ] \right  \} \nonumber \\  
 L_{0} & = & 2 N_{1} +  N_{2} - N_{4}  - 2 N_{5} =  2 H_{1} +  
 \frac{3}{2} H_{2}
 \eea 
 where $H_{1}, H_{2}$ of  Cartan generators of $so_{q}(5) $ and
 we have used eqs.(\ref{eq:real})-(\ref{eq:song}) and the deforming map eq.(\ref{eq:CZ}).
 
  \subsection{ \bf $q$-analogue of the  rotational embedding chain} 

\bea
gl_{q}(6) \supset gl_{q}(3) & \nonumber \\
& \supset^{\Ua_{q \ra 1}}  so_{q}(3)  
\label{eq:dr}
\eea
the $q$-boson realization of  $sl_{q}(3) \subset gl_{q}(6)$, \cite{NS}
 \bea
E_{1}^{\dag}  &  =  & \left \{q^{N_{4} - N_{2}/2}  \, \sqrt{q^{N_{1}} + q^{-N_{1}}}Ê\,  
b_{1}^+ b_{2}  +  b_2^+ b_4 \, q^{N_{1} - N_{2}/2} 
\sqrt{q^{N_{4}} + q^{-N_{4}}}Ê \right \} \, q^{-(N_{3} - N_{6})/2}  
\nonumber \\Ê 
&  +  &  b_3^+ b_5 \, q^{(N_{1} - N_{4})}  \nonumber \\Ê
E_{2}^{\dag}  &  =  & \left \{q^{N_{6} - N_{5}/2} \,\sqrt{q^{N_{4}} + q^{-N_{4}}}Ê\,  
b_4^+ b_5 \, +  \, b_5^+ b_6 \, q^{N_{4} - N_{5}/2} 
\sqrt{q^{N_{6}} + q^{-N_{6}}}Ê \right \} \, q^{(N_{2} - N_{3})/2} \nonumber \\Ê 
& +  & b_2^+ b_3 \, q^{-(N_{4} - N_{6})}  \nonumber \\  
 H_{1} & = & 2 N_{1}  - 2 N_{4} + N_{3}  - N_{5} \nonumber \\ 
  H_{2} & = &  N_{2}  -  N_{3} +  2 N_{4}  - 2 N_{6}    
\eea
the generators of $so_{q}(3)$ can be written in terms of the 
generators of $gl_{q}(6)$  
\bea
& L_{+}  =  {\cal Q_{+}}  \left \{ 2 \left [ 
\sqrt{\frac{N_{1}}{[N_{1}]_{q}}} e^{+}_{1} 
\sqrt{\frac{N_{2}}{[N_{2}]_{q}}}  
+ \sqrt{\frac{N_{2}}{[N_{2}]_{q}}} \, [e^{+}_{2}, e^{+}_{3}]_{q} \, 
\sqrt{\frac{N_{4}}{[N_{4}]_{q}}} \, q^{N_{3}} \, \right ]  \right. 
\nonumber \\  
& + \sqrt{2} \left [  \sqrt{\frac{N_{3}}{[N_{3}]_{q}}} \, [e^{+}_{3}, e^{+}_{4}]_{q} 
\, \sqrt{\frac{N_{5}}{[N_{5}]_{q}}} q^{N_{4}}  + \sqrt{\frac{N_{2}}{[N_{2}]_{q}}} \, 
e^{+}_{2} \, \sqrt{\frac{N_{3}}{[N_{3}]_{q}}} \, \right ] \nonumber \\  
& +  2 \, \left. \left [ \sqrt{\frac{N_{4}}{[N_{4}]_{q}}} e^{+}_{4} \sqrt{\frac 
{N_{5}}{[N_{5}]_{q}}} \, + \, \sqrt{\frac{N_{5}}{[N_{5}]_{q}}} e^{+}_{5} 
\sqrt{\frac{N_{6}}{[N_{6}]_{q}}} \right]  \, \right \} \nonumber \\ 
& L_{0}  =   2 N_{1} +  N_{2} - N_{5}  - 2 N_{6} =  H_{1} +   H_{2}
 \eea
 where $H_{1}, H_{2}$ of  Cartan generators of $sl_{q}(3) $

 \subsection{ \bf $q$-analogue of the $\gamma$-unstable embedding chain } 

\bea
gl_{q}(6) \supset so_{q}(6) & \nonumber \\
 & \supset^{\Ua_{q \ra 1}} so_{q}(5)  & \nonumber \\
 &  \mbox{}  &    \supset^{\Ua_{q \ra 1}} so_{q}(3)  
\label{eq:dg}
\eea
the $q$-boson realization of  $so_{q}(6) \subset sl_{q}(6)$ \cite{NS}
 \bea
E_{1}^{\dag} &  =  & b_2^+ b_{4} \, q^{(N_{3} - N_{5})/2} \,  
 +  \,  b_3^+ b_{5} \, q^{-(N_{2} - N_{4})/2}   \nonumber \\Ê 
E_{2}^{\dag} &  =  & b_1^+ b_{2} \, q^{(N_{5} - N_{6})/2} \,  
 +  \,  b_5^+ b_{6} \, q^{-(N_{1} - N_{2})/2}    \nonumber  \\Ê 
 E_{3}^{\dag} &  =  & b_2^+ b_{3} \, q^{(N_{4} - N_{5})/2} \,  
 +  \,  b_4^+ b_{5} \, q^{-(N_{2} - N_{3})/2} \nonumber \\  
 H_{1} & = & N_{2}  - N_{4} + N_{3}  - N_{5}  \nonumber \\  
  H_{2} & = & N_{1}  - N_{2} + N_{5}  - N_{6} \nonumber \\  
   H_{3} & = & N_{2}  - N_{3} + N_{4}  - N_{5}  
\eea
the $q$-boson realization of  $so_{q}(5)$, deformation of $SO(5)$ maximal 
subalgebra of $SO(6) \subset SU(6)$, using eq.(\ref{eq:map}), can be written
\bea
& E_{1}^{\dag}  =  \left (\sqrt{\frac{N_{1}}{[N_{1}]_{q}}} e^{+}_{1} 
\sqrt{\frac{N_{2}}{[N_{2}]_{q}}}  +  \sqrt{\frac{N_{5}}{[N_{5}]_{q}}} e^{+}_{5} 
\sqrt{\frac{N_{6}}{[N_{6}]_{q}}} \right )
\sqrt{\frac{[\widehat{H}_{1} + \widehat{H}_{2}  +  1]_{q} \, [\widehat{H}_{2}]_{q}}
{(\widehat{H}_{1}  + \widehat{H}_{2}  +  1)  \widehat{H}_{2}}}
\nonumber \\
& E_{2}^{\dag}  =  \left (\sqrt{\frac{N_{2}}{[N_{2}]_{q}}} \, [e^{+}_{2}, e^{+}_{3}]_{q} \, 
\sqrt{\frac{N_{4}}{[N_{4}]_{q}}} \, q^{N_{3}} + 
\sqrt{\frac{N_{3}}{[N_{3}]_{q}}} \, [e^{+}_{3}, e^{+}_{4}]_{q} 
\, \sqrt{\frac{N_{5}}{[N_{5}]_{q}}} q^{N_{4}}  \right. \nonumber \\
& + \left. \sqrt{\frac{N_{2}}{[N_{2}]_{q}}} \,e^{+}_{2} \, \sqrt{\frac{N_{3}}{[N_{3}]_{q}}}
+ \sqrt{\frac{N_{4}}{[N_{4}]_{q}}} e^{+}_{4} \sqrt{\frac 
{N_{5}}{[N_{5}]_{q}}} \right )
\nonumber \\
& \times   \frac{2}{q + q ^{-1}} \, 
\sqrt{\frac{[\widehat{H}_{2}  +  1]_{q} [- \widehat{H}_{2}  -  2)]_{q}}
{( \widehat{H}_{2}  +  1) (- \widehat{H}_{2}  -  2)}}
\eea
\be
 \widehat{H}_{1} = N_{1} + N_{5} - N_{2} - N_{6} \;\;\;\;\;\; 
 \widehat{H}_{2} = 2 (N_{2} - N_{5}) 
\ee
where $\widehat{H}_{1}, \;\; \widehat{H}_{2}$ are the Cartan generators of $so_{q}(5)$ 
and the generators of $so_{q}(3)$ can be written in terms of the 
generators of $gl_{q}(6)$
\bea
& L_{+}  =  {\cal Q_{+}} \left \{ \sqrt{2}\left [ 
\sqrt{\frac{N_{1}}{[N_{1}]_{q}}} e^{+}_{1} \sqrt{\frac{N_{2}}{[N_{2}]_{q}}}
  +   \sqrt{\frac{N_{5}}{[N_{5}]_{q}}} \, e^{+}_{5}  \, 
\sqrt{\frac{N_{6}}{[N_{6}]_{q}}}  \, \right ]  \right. \nonumber \\  
& + \left.   \sqrt{\frac{3}{2}} \, \left [\sqrt{\frac{N_{2}}{[N_{2}]_{q}}} \,
[e^{+}_{2}, e^{+}_{3}]_{q} 
\, \sqrt{\frac{N_{4}}{[N_{4}]_{q}}} q^{N_{3}} +
\sqrt{\frac{N_{3}}{[N_{3}]_{q}}} \,[e^{+}_{3}, e^{+}_{4}]_{q} \, 
\sqrt{\frac{N_{5}}{[N_{5}]_{q}}} q^{N_{4}} \right. \right. \nonumber \\  
& +  \left. \left.  \sqrt{\frac{N_{2}}{[N_{2}]_{q}}} \, e^{+}_{2}  \, 
\sqrt{\frac{N_{3}}{[N_{3}]_{q}}} + \sqrt{\frac{N_{4}}{[N_{4}]_{q}}} \, e^{+}_{4}  \, 
\sqrt{\frac{N_{5}}{[N_{5}]_{q}}} \right ] \right \}\nonumber \\  
& L_{0}  =   2 N_{1} +  N_{2} - N_{5}  - 2 N_{6} = \frac{3}{2} ( 
 H_{1}  + H_{3}) +  2 H_{2}
 \eea  
where $H_{1}, H_{2},  H_{3}$ of  Cartan generators of $so_{q}(6) $. We 
recall that
\be
 \widehat{H}_{1} = H_{2}   \;\;\;\;\;\;\;\;\;\;\;\; 
 \widehat{H}_{2} =  H_{1} - H_{3} 
\ee

\section {\bf   Conclusions}

Starting from a spectrum generating algebra
\be
U(6) \supset L \supset SO(3)
\ee
We have shown that, using $q$-bosons, a deformed analogue of this  chain can 
be obtained  replacing $U(6)$ by $gl_{q}(6)$ and writing the 
generators of $L_{q}$ ($L = SO(6), SO(5),$ 

$SU(3), SO(5) \subset SO(6) $) and of $so_{q}(3)$ in terms of the generators 
of $gl_{q}(6)$  (or of  $gl_{q}(5)$  for the vibrational chain)
with coefficients taking value in $gl_{q}(6)$ (or $gl_{q}(5)$). In order to write our 
results we have used the relation between standard bosonic operators and 
 $q$-bosons eq.(\ref{eq:song}), the deforming  map between $su(2)$ and  $su_{q}(2)$
 eq.(\ref{eq:CZ}) and 
   between $sp(4) \equiv so(5)$ and  $sp_{q}(4) \equiv so_{q}(5)$ 
   eq.(\ref{eq:map}).
Remark that in the last step of the embedding chain one can keep   undeformed $SO(3)$ and the 
generators $l_{\pm}$ of the undeformed algebra of the angular momentum can be written, 
e.g., as 
\be
 l_{+} = \sum_{k} \, A_{k}^{+}(q, 
\{e_{i}^{+},h_{i}\}) \, e_{k}^{+}   
\ee 
where $A_{k}^{+}(q, \{e_{i}^{+},h_{i}\})$ are operators taking value 
in $gl_{q}(6)$ (or $gl_{q}(5)$).  
Let us point out that the problem in a deformation of the subalgebra
of a subalgebra in an embedding chain of the type, e.g.,  
 \be
gl_{q}(6) \supset  L_{q} \supset  so_{q}(3)  
\ee 
is that the generators $L_{\pm}$ cannot be written in terms of the 
generators of $L_{q}$, but only in terms of the generators of 
$gl_{q}(6)$. This peculiar feature is common to any deformed algebra 
when we consider embedding chain with deformed singular subalgebra
$ L \subset G$.
In conclusion we have shown that it is possible to write the Cartan-Chevalley 
generators
of $so_q(6)$ and $sl_q(3)$ in terms of the generators of $gl_q(6)$ and those
of $so_q(5)$ in terms of $gl_q(5)$, but that it is not possible to extend
further the procedure to obtain, in particular, $so_q(3)$, as it was argued in  
\cite{Scia2} and \cite{NS}. However for the considerd chains 
eq.(\ref{eq:AI}) we can go a step further writing the generators in 
terms of those of the {\em grand-mother} $gl_{q}(6)$. Even if we have considered a 
particular case, so our results are not quite general, the used procedure is enough 
general to  be applied successfully to other physically relevant models, even  
  taking into account supersymmetric extensions.
In the  spirit of the use of spectrum generating algebra, one should 
write a hamiltonian of a physical system as a sum of invariants of the
$q$-algebras appearing in the embedding chain of the previous Section.
This can be done as the Casimir operators of $sl_q(n)$  \cite{Cha1},
\cite{Bin} and of $so_q(5)$ \cite{ZHB}, \cite{Cha2} are known. For
physical application one needs self-adjoint Casimir. For $q$ real this
property is guaranteed in the realization we have used. In many 
applications of $q$-algebras, however, it seems that $q$ being a phase
is a favoured value, see \cite{BD}. For this value, due to the lack of
invariance $ q \; \lra \; q^{-1} $ in our expressions,   the product
$E^+ \, E^-$ is no more self-adjoint. In order to restore the 
self-adjointness of the hamiltonian one has to sum the Casimir operators of  $L_{q}$ 
($so_{q}(3)$)  and of $L_{q^{-1}}$  ($so_{q^{-1}}(3)$). As a final
comment, if one considers transitions in the physical system induced by
$q$-tensor operator under $so_{q}(3) $ or, if the rotation group in 
physical space is left undeformed, by tensor operator under $SO(3) $,
one has to face the problem of the choice of the coproduct. In order 
to have the usual structure of $q$-tensor eq.(\ref{eq:tensor}) or 
 tensor  operators,  eq.(\ref{eq:tensor})  in the limit $q \ra 1$,
the coproduct $\Delta$  has to be imposed in $ so_{q}(3)$ or $SO(3)$  
 and cannot be inherited by that  of   $G_{q}$ or $L_{q}$ .
    
\bigskip

 \textbf{Acknowledgments}:  This work is an extended version of the talk 
 given in the XII Int. Symposium: {\em Symmetries in Science}, Bregenz 
 22-27 July 2001. I thank the organizers Prof. B. Gruber and Prof. M. 
 Ramek for their invitation. This work  has been partially supported by the  
 M.U.R.S.T.  through  National Research Project {\sl SINTESI 2000}.

\end{document}